Suppression of random noise by the separation of frequencies

V. N. Tibabishev

(asvt51@narod.ru, tivasl@vaush.net)

**It is shown that a well-known theory of random stationary processes contain contradictions. Integral representations of correlation functions and random stationary processes are investigated further. The new method of struggle with handicaps is received on the basis of the carried out researches. Method of dealing with noise leads to a new method of identification of dynamic characteristics of control objects in a class of multidimensional linear stationary models. As an example, describe an algorithm for obtaining the differential equation-wire feed control aircraft pitch, taking into account the elastic deformation of the structure of class IL - 96 aircraft.**

## 1. Introduction

The initial (primary) data are given in the presence of noise in the solution of applied problems always. Certain properties of the signals used to control noise. The modern theory of random stationary in the wide sense with finite variance processes (the random process) is described in many books, for example, in [1, 2, 3and 4]. In the theory of random processes used empirical [5, pp. 567-623] and theoretical methods, for example, published in [6] [7], [2]. Empirical studies of E. Slutsky summarized, for example, in [4, pp. 337, 340]. From the results of empirical studies suggest that a stationary random process and its correlation function is represented by a trigonometric Fourier series of arbitrary frequency, and therefore belong to the Hilbert space of almost periodic functions of H. Bohr $B_2(-\infty,+\infty)$. In [2, p. 166] indicated that the random processes with discrete spectrum does not possess the ergodic property.

The first theoretical concept was published by A. Khinchin, according to which the correlation function of random process has two components [6, p. 49]. One component belongs to the Hilbert space of almost periodic functions $B_2(-\infty,+\infty)$. Average value of the time of the square the other component is always zero. This implies that another component belongs to the Hilbert space of square-integrable Lebesgue functions $L_2(-\infty,+\infty)$.

Another theoretical concept was published by Wiener, according to which the correlation function of a stationary process belongs only to the Hilbert space of square integrable functions, Lebesgue integral is the Fourier-Plancherel terms of the spectral density [7, p. 275]. Such a model is a random process has the ergodic property of the theorem, published in [2, c. 164].

In the works of Norbert Wiener and Khintchine representation of random processes themselves were not considered. This problem was solved empirically first E. Slutsky, then in theory by Kolmogorov [8, p. 340], further Cramer [2, p. 135] and other researchers. In all views, except of E. Slutsky, a stationary random process is an integral of the Fourier-Stieltjes.

All three concepts of a random stationary process contradict each other. In Wiener [9, p. 193] shows that if the process belongs to the Hilbert space $L_2(-\infty,+\infty)$, then the correlation function, obtained by averaging over time is identically zero on the entire line. This assertion N. Wiener is contradiction ergodic theorem given in [2, c. 164].

Middleton, D. [3, p. 191] and other researchers [10, pp. 88] pointed out that the integral of the square of a random stationary process for all direct costs, even when the average value of zero. Therefore, a stationary random process can not belong to a Hilbert space $L_2(-\infty,+\infty)$.

Known [11, c. 375], that the space $L_2(-\infty,+\infty)$ is a set of functions that decay sufficiently fast at infinity. Such functions are not suitable for description of random stationary processes on the whole line. In many monographs as examples of stationary processes are functions that do not belong to the space $L_2(-\infty,+\infty)$, a are elements of almost periodic functions $B_2(-\infty,+\infty)$ [12, pp. 50], [2, pp. 166], [3, pp. 60], [4, pp. 336]. It is still unclear why the Lebesgue-Fourier-Stieltjes in N. Wiener [7, pp. 274] has led to conflicting results.

A hypothesis, that the contradictions are due to incorrect assumptions adopted in theoretical studies. First, all theoretical studies suggest that the angular frequency $\omega$ is a continuous variable in the integration [8, pp. 243] and differentiation [7, c. 274]. Second, the probability distribution function contains a step and continuous components [7, c. 274], [2, c. 25], [6, c. 49]. In this connection there arose the first task of assessing the correctness of these assumptions.

## 2. Investigation of the known assumptions

First, we show that a certain assumption about the continuity of the independent variable angular frequency is incorrect assumption. Angular frequency $\omega$ is related to the cyclic frequency $f$ numerical function $\omega = 2\pi f$. Cyclic frequency belongs to a continuous set of real numbers $f \in W$. The answer to the question, to what set of numbers belongs to the set of corner frequencies, the following theorem.

*It is theorem 1. Let the real function $y = kx$ defined on the set of real numbers $x \in W$ for all $x \neq 0$. If the ratio $k$ is an element of a subset of irrational numbers, transcendental $k \in T$, for example, $k = \pi$ the numerical function that displays the specified set of real numbers $x \in W = V \bigcup I \bigcup T$, only a subset of irrational numbers, transcendental $y \in T$, where $V$ - a subset of rational numbers $I$ - a subset of algebraic irrationals, for example, $\sqrt{2}$ where the union $V \bigcup I$ - is an infinite countable set of algebraic numbers.*

It is proof. Known [13, pp. 17], that there exists $x \neq 0$ an inverse $x^{-1}$ for each element such that $x^{-1}x = 1$. Multiply the left and right sides of a numerical function to return the item to an arbitrarily chosen $x \neq 0$. Obtain $x^{-1}y = kx^{-1}x \equiv k \in T$. It follows that for $x^{-1} \neq 0$ an arbitrary product $x^{-1}y \in T$ of an irrational transcendental number. Of all the possible random numbers $x^{-1}$ choose an arbitrary rational number $x^{-1} = m/n$, where m and n are arbitrary nonzero integers. Of all the random numbers $m$, choose a number $m = n$. In this case, we find that $y = k \in T$. Theorem 1 is proof.

The question arises, under what factors $k$ numerical function $y = kx$ defined on the real line $x \in (-\infty,+\infty)$ for all $x \neq 0$ will display the specified set of real numbers into itself ? The answer to this question is the following theorem.

*It is theorem 2. Let the real function $y = kx$, defined on the set of real numbers $x \in W$ for all $x \neq 0$. If the coefficient $k \neq 0$ belongs to a subset of rational numbers $k \in V \subset W$, then the numerical function displays the specified set of real numbers $x \in W$ into itself $y \in W$.*

It is proof. All rational coefficients $k \neq 0$ are presented as a fraction $k = m/n$, where $m$ and $n$ are arbitrary nonzero integers. Since the choice take $m = n$. In this case $y = x \in W$, as required.

It is consequence of theorem 2. Angular frequency $\omega = \pi 2 f$, where $f \in W$ - angular frequency, which belongs to the set of real numbers. It is according with theorem 2 the set of numbers $2f \in W$. Since $\pi \in T$ it belongs to the set of irrational numbers, transcendental, then by theorem 1 the set of corner frequencies $\omega \in T$ belong to the subset of irrational transcendental number, which is not continuous set of numbers. In an infinite uncountable interval only irrational transcendental numbers contains an infinite countable number of seats are not engaged in irrational, transcendental numbers. On an infinite line of real numbers $-\infty < \omega < +\infty$, these seats are a countably infinite set of algebraic numbers [13, pp. 48]. It follows that an infinite range of irrational transcendental frequency contains an infinite countable set of discontinuity points. On an infinite line of real numbers, these seats are a countably infinite set of algebraic numbers [13, pp. 48]. It follows that an infinite range of irrational transcendental frequency contains an infinite countable set of discontinuity points. On a subset of irrational numbers, transcendental notions $\omega \in T$ of mathematical analysis as an infinitely small quantity $d\omega$, the derivative (spectral density) and the differential of the independent variable angular frequency $\omega$ does not exist. With respect $\omega$ to the integral sums of Riemann, Riemann-Stieltjes and Fourier and Fourier-Plancherel not exist.

Second, the well-known assumption that the distribution function of angular frequency contains a continuous component, is an invalid assumption.

It is theorem 3. *In general, the probability distribution as a function of real random variables consists of three components - a step function, absolutely continuous differentiable function and a singular function. If the probability distribution is a function of random angular frequency, defined on the set of irrational numbers, transcendental, then the probability distribution function of random angular frequency contains only a step function with an infinite countable number of discontinuity points and contains no nonzero continuous components.*

It is proof. It is known [2, pp. 25], that in general the distribution function of random real variables $x \in (-\infty, +\infty)$ consists of three components

$$F(x) = F_1(x) + F_2(x) + F_3(x),$$

where $F_1(x) = F_1(x_k), x_k \le x < x_{k+1}$ - a step function of discrete values $x_k$, $k = 0,1,2,3,...$, continuous variable $x$, and $x_k$ can also take an infinite value, $F_2(x)$ - absolutely continuous differentiable function on a continuous variable $x$, $F_3(x)$ - a singular function - a continuous function, which has a derivative almost everywhere equal to zero, for example, $F_3(x) = c$ where $c$ - arbitrary constant.

From the properties of the probability distribution function defined, for example, random variables $x \in (-\infty, +\infty)$, it follows that the entire function $F(-\infty) = 0$ and its two components $F_1(-\infty) = 0$ and $F_2(-\infty) = 0$. At the same time $F_1(-\infty) + F_2(-\infty) + c = 0$. It follows that $c = 0$. The distribution function is represented as $F(x) = F_1(x) + F_2(x)$. From the other properties of the probability distribution of the variable $x$ that $F(\infty) = F_1(\infty) + F_2(\infty) = 1$. Hence

$$F_1(\infty) = 1 - F_2(\infty) . \quad (1)$$

It is known [14, pp. 30], that the step function can have a finite number of steps, for example, $n$ or an infinite number of steps, but necessarily countable. Therefore, the component $F_1(\infty)$ belongs to a countable set of rational numbers and can take values in the interval [0,1]. Let $F_1(\infty) = n/m$ where $n$ and $m$ - integers with $n \le m$ .

Into an infinite number of steps tends to its maximum possible value of the countable $n \to m = \infty$ . In this case, a discrete component tends to its maximum value $F_1(\infty) = 1$ and the continuous component tends to its minimum value $F_2(\infty) = 0$ . Since when the number $\omega = \infty$ of breaks in the step function equal to an infinite countable set of values, for all $\omega$ the distribution function continuous component $F_2(\omega) \equiv 0$ . Theorem 3 is proved.

Thus, the probability distribution function of angular frequency has only a discrete component, containing an infinite but countable set of stairs. In this case, it contains no nonzero continuous components.

### 3. Representation of correlation functions and processes

In this section we repeat the well-known findings of Wiener [7, pp. 274] to represent the correlation functions and the known output representation for a random process, Cramer [2, pp. 138] with the difference that we use is not anticipated by the authors properties for the angular frequency and its distribution function, as proved above properties.

Corollary 1 is of theorem 3. If the distribution function of the random variable contains a step function with an infinite countable number of steps, then the Lebesgue measure, generated by this distribution function contains only discrete measure. Therefore, the components of the Lebesgue-Stieltjes measure on a continuous non-existent.

Since the angular frequency $\omega$ is not a continuous variable, the function $\exp(j\omega t)$ of the variable $\omega$ is not continuous. Fourier-Stieltjes integral is defined only for continuous functions [14, c. 82]. Since this condition is not fulfilled, instead of the Riemann-Stieltjes integral must be taken Lebesgue-Stieltjes integral [11, pp. 337], as did Norbert Wiener

$$R(t) = \int_{-\infty}^{+\infty} \exp(j\omega t) d\mu_F = \int_{-\infty}^{+\infty} \exp(j\omega t) dF(\omega) \ , \quad (2)$$

where $\mu_F$ - the measure generated by the distribution of the random corner frequency $F(\omega)$ on a subset of transcendental numbers $\omega \in T$ .

Known [11, pp. 337], that in cases with a discrete measure $\mu_F$ of the Lebesgue-Stieltjes integral (6) leads to nonintegral concepts defined in the space of almost periodic functions

$$R_1(t) = \sum_{k=-\infty}^{k=+\infty} \exp(j\omega_k t) h(\omega_k) , \quad (3)$$

where $h(\omega_k)$ - at a frequency jump $F(\omega)$ function $\omega_k$ .

Submission of a random process differs from the representation of the correlation function so that the integral representation (2) instead of the distribution function $F(\omega)$ uses a random complex continuous with bounded variation with zero mean and orthogonal increments spectral process $\zeta(\omega)$ [2, pp. 138]

$$\xi(t) = \int_{-\infty}^{+\infty} \exp(j\omega t) d\zeta(\omega) . \qquad (4)$$

It is corollary 2 of Theorem 3. It is known [14, pp. 41], that the function of bounded variation in a certain way is the sum of functions of the jumps and the amount of continuous distribution functions. Therefore Theorem 3 is valid not only for the non-decreasing function of probability apply, but also for functions with limited modifications.

It is corollary 3 of Theorem 3. Since any function of bounded variation is the difference between two step and the difference of two continuous distribution functions [14, c. 41], the function of bounded variation of the angular frequency $\zeta(\omega)$, which does not contain continuous components, can generate only a discrete measure $\mu_F = d\zeta(\omega)$.

Function $\exp(j\omega t)$ of the variable $\omega$ contains an infinite countable set of discontinuity points. Therefore, for integration into the mapping (4) must be taken Lebesgue-Stieltjes integral. By corollary 3, Embedding theorems 3 Spectral process $\zeta(\omega)$, which contains an infinite countable set of discontinuity points, can produce only a discrete measure $\mu_F = d\zeta(\omega)$. Known [11, pp. 338] that at least a discrete integral representation (4) becomes complex Fourier series, defined in the Hilbert space of almost periodic functions of Bohr $B_2(-\infty,+\infty)$

$$\xi(t) = \sum_{k=-\infty}^{k=+\infty} C(j\omega_k) \exp(j\omega_k t) , \qquad (5)$$

where $\omega_k$ - the point of tears.

Our representation for correlation functions (3) and for stationary random processes (5) are differ from the known representations received by Wiener and G. Kramer. First, a random process and its correlation function does not belong to Hilbert space $L_2(-\infty,+\infty)$ and the Hilbert space of almost periodic functions $B_2(-\infty,+\infty)$. In this space integral of the square of the random process does not diverge. Second, a stationary random process and its correlation function has only discrete spectrum. In this theoretical study are consistent with empirical studies E. Slutsky.

## 4. Functional space for random stationary processes

We shall define a functional space for random stationary processes with discrete spectra on the example of the object control, containing the $l$ inputs and $d$ outputs. It is believed that a multi-dimensional control object allows $l \cdot d$ an approximate description of the dynamic performance of control channels in a class of linear stationary models. Original raw data presented in the form of single realizations of random nonergodic processes only with discrete spectra, which simultaneously observed (recorded) on the inputs $\tilde{x}_{iq}(t)$ and output $\tilde{y}_{ip}(t)$ of control object, where $i = 1,2,3,...$ - number of implementation $q \in [1,l]$ - the number of input $p \in [1,d]$ - output number. Instead of accurate baseline data $x_{iq}(t)$ and $y_{ip}(t)$ can be observed or synchronous writes only approximate data, $\tilde{x}_{iq}(t) = x_{iq}(t) + n_{iq}(t)$ and $\tilde{y}_{ip}(t) = y_{ip}(t) + m_{ip}(t)$ in which exact components distorted by random noise or additive $n_{iq}(t) \in N_q$ and $m_{ip}(t) \in M_p$.

Each implementation of approximate initial data belongs to the Hilbert space of almost periodic functions. There are various generalizations of almost periodic functions [15, pp. 218]. We will consider such control objects, whose dynamics obeys the laws of classical mechanics, in particular, second order differential equation of Newton, for example, in the form of $c\ddot{y}(t) = f(t)$

, $\ddot{y}(t)$ - the second derivative, for example, the output signal , $f(t)$ - force. According to the known existence theorem [16, pp. 150] function $\ddot{y}(t), f(t)$ must be continuous. In many cases [16, pp. 233] the function $f(t)$ is the Fourier series. Since the functions $\ddot{y}(t), f(t)$ are continuous functions, the Fourier series for these functions are convergent series in the known space [17, pp. 341]. In the Hilbert space of almost periodic functions in the sense of Besicovitch each element appears to be convergent Fourier series [15, pp. 223]. Therefore, we assume that each component of the approximate initial data are elements of the Hilbert space of almost periodic functions in the sense of Besicovitch $B_2(-\infty,+\infty)$.

We consider the ergodic properties of random processes in this space. Let the stationary random process $x(t)$ with bounded non-zero expectation $E\{x(t)\} = m$ belongs to the Hilbert space of almost periodic functions in the sense of Besicovitch $x(t) \in B_2(-\infty, +\infty)$. It is known [12, pp. 184], that in this space, the notion of average

$$M\{x_i(t)\} \equiv \lim_{T\to\infty} \frac{1}{2T} \int_{-T}^{T} x_i(t)dt,$$

where $i$ - number of implementation.

Known ergodic theorem Birkhoff - Khinchin proved that for non-Hilbert space of integrable functions on Riemann [2, pp. 158], that is, for the space $L_1(-\infty, +\infty)$. This theorem is generalized to the Hilbert space of almost periodic functions in the sense of Besicovitch.

It is theorem 4. *Suppose that, in general, a stationary random process belongs to the set of random processes $x_i(t) \in X_m$ with limited non-zero expectation $E\{X_m\} = m$. If a random stationary process is an element of the Hilbert space of almost periodic functions in the sense of Besicovitch, then with probability 1, the average value of any implementation of this set $X_m$ equal to the mean stationary random process.*

It is proof. Every realization of a random process from a variety of $x_i(t) \in X_m = \{x_i(t) \in X_m\colon E\{X_m\} = m\}$ can be represented as $x_i(t) = m + \check{x}_i(t)$ , where $\check{x}_i(t)$ - centered stationary random process. Each centered stationary random process that belongs to the Hilbert space of almost periodic functions in the sense of Besicovitch is convergent trigonometric series [4, pp. 340]

$$\check{x}_i(t) = \sum_{k=1}^{k=\infty}(a_i(\omega_k)cos(\omega_k t) + b_i(\omega_k)\sin(\omega_k t)), \tag{6}$$

where $\omega_k > 0$.

It is $M\{x_i(t)\} = E\{X_m\} = m$ since with probability 1, the average value of a centered random process is zero. It follows that a stationary random process with discrete spectrum has the ergodicity of the first order [1, pp. 133]. Theorem 4 is proved.

We show that the second and mixed moments, e.g., centered random stationary processes belonging to the Hilbert space of almost periodic functions in the sense of Besicovitch, Ergodic properties of second order do not possess. Let $\check{x}_i(t) \in X_\rho$, where $X_\rho$ - the set of almost periodic functions by generating a correlation function of a certain species [4, pp. 336]

$$R_x(\tau) = \sum_{k=1}^{k=\infty} E\{C_i^2(\omega_k)\}\cos(\omega_k \tau) , \tag{7}$$

where $E\{C_i^2(\omega_k)\} = E\{a_i^2(\omega_k)\} + E\{b_i^2(\omega_k)\}$ .

The operator display the space $B_2(-\infty, +\infty)$ into itself [12, pp. 186], because it is completely continuous operator of convolution type [12, pp. 186]. Let the operator equation of convolution type $Ax_i = z_i$ kernel is generated by a function of $x_i(t-\tau)$ . In this case, we find that

$$z_i(\tau) = M\{x_i(t-\tau)x_i(t)\} = \sum_{k=1}^{k=\infty} C_i^2\,(\omega_k)\cos(\omega_k\tau), \quad (8)$$

This shows that the stationary processes that belong to the Hilbert space of almost periodic functions in the sense of Besicovitch, do not possess the second-order ergodicity as $R_x(\tau) \neq z_i(\tau)$.

The exact component of the input signal generates a forced movement control channels. Noise, distorting the exact component of the input signal, no effect on the forced movement of the control channels. The output signals induced motion is distorted by other noise not associated with precise input to any obstacle that distorts the exact input. Interference and accurate random processes are of different nature. If the components are the exact image and inverse image of the observed signals of linear time-invariant operators, then the interference is linearly independent and therefore uncorrelated processes. In the theory of random processes produce correlated and uncorrelated random processes. Uncorrelated random processes in the sense of the conditions $E\{x_{iq}n_{iq}\} = 0$ , $E\{y_{ip}n_{iq}\} = 0$ , $E\{y_{ip}m_{ip}\} = 0$ , $E\{x_{iq}m_{ip}\} = 0$ , $E\{m_{ip}n_{iq}\} = 0$ , for all $i,q$ and $p$ , where $E\{.\}$ - the symbol of averaging over infinite-dimensional set pieces implementations. If, for example, the control channel between the input and output q p is the control channel, the exact components have a nonzero cross-correlation $E\{x_{ip}y_{iq}\} \neq 0$ . In the theory of random processes, these differences in the correlations used to combat interference. To check that these conditions are necessary baseline data in the form of sets containing an infinite number of realizations. It is practically impossible.

We consider the conditions under which, firstly, the primary source data given in the form only isolated implementations and, secondly, the single realization of random processes do not possess the ergodic property. Receipt of secondary source data in the form of deterministic correlation functions in such circumstances is impossible. Therefore, we use the condition of linear independence of the individual realizations, under which condition is uncorrelated to be checked on a set of implementations, is always satisfied.

Consider the equation for the submission of a centered random process (6). In this expression, the amplitude of random $a_i(\omega_k)$ and $b_i(\omega_k)$ deterministic harmonic oscillations with frequencies $\omega_k$ generate lots of random processes. From this expression it follows that all of a given set of random processes is conceived as a whole because it is the linear hull spanned by a deterministic basis of harmonic functions whose frequencies coincide with frequencies of the harmonics of the correlation function $R_x(\tau)$. The set of frequencies of harmonic components of the observed signals is a deterministic secondary source data. The main difference between the proposed method for dealing with noise from the known method is that as a secondary information are not used correlation functions, and set the frequency of the harmonic components. It turned out that the properties of systems of sets of frequency harmonic components of the observed signals depend on the correlation (linear dependence) constitute the observed signals.

## 5. Properties of independent random processes

Distinguish multidimensional control objects with linearly dependent (correlated) and linearly independent (uncorrelated) input actions. First, consider the control objects with linearly independent input actions. You can specify a number of linearly independent systems of sets of random processes. In particular, a system of linearly independent sets of exact input actions $S_x = \bigcup_{q=1}^{q=l} X_q$ where $x_{iq}(t) \in X_q$ - exact sets of random input processes (impacts). Elements $x_{iq} \in X_q$ are linearly independent if the equality $\alpha x_{i1} + \beta x_{i2} + \cdots + \gamma x_{il} = 0$ follows from the equality $\alpha = \beta = \cdots = \gamma = 0$ for all $i$ .

System of linearly independent sets of exact input signals and additive noise, distorting the accurate baseline data $S_{xnm} = \bigcup_{q=1}^{q=l} X_q \bigcup_{c=1}^{c=l} N_c \bigcup_{b=1}^{b=d} M_b$ , where $n_{ic}(t) \in N_c$, $m_{ib}(t) \in M_b$. Each set of a lot of thought as a whole due to the fact that $S_{xnm}$ averaging over the set of products centered random processes, such as

$$x_{iq}(t) \in X_q = \{x_{iq} : E\{x_{iq}(t-\tau)x_{iq}(t)\} = R_{xq}(\tau)\}$$

generates a deterministic correlation function $R_{xq}(\tau)$.

In systems of sets, for example, $S_{xnm}$ introduce a scalar product for complex centered random functions $x_{iq}(t) \in X_q$ and $x_{ip}(t) \in X_p$, for example,

$$(x_{iq}(t), \bar{x}_{ip}(t)) \equiv \lim_{T\to\infty} \frac{1}{2T} \int_{-T}^{+T} x_{iq}(t)\bar{x}_{ip}(t)dt = M\{x_{iq}(t)\bar{x}_{ip}(t)\}$$

where $x_{ip}(t)$ and $\bar{x}_{ip}(t)$ - complex conjugate functions.

It is known [12, pp. 185] that uncountable harmonic basis has a non-separable complex Hilbert space of almost periodic functions $B_2(-\infty,+\infty)$. In this space, each vector is only a countable sum of nonzero orthogonal projections. Therefore, we assume that, for example, $S_x \subset B_2(-\infty,+\infty)$ is a separable subspace of complex functions. In this regard, given the above mathematical expression is translated into the complex plane.

By hypothesis, any two simultaneously recorded sale $x_{iq}(t) \in X_q$ and $x_{ip}(t) \in X_p$ must be linearly independent for all $q \neq p$. It is known [11, p.136], if the vectors $x_{iq}(t) \in X_q$ and $x_{ip}(t) \in X_p$ are orthogonal, then they are linearly independent. This implies another condition of linear independence of random realizations of random processes ( $x_{iq}(t), \bar{x}_{ip}(t)) = 0$ , where $x_{ip}(t), \bar{x}_{ip}(t)$ - complex conjugate functions.

From formula (6) shows that for each $i$ - th realization of a random process, belonging to a particular set of random processes $x_{iq}(t) \in X_q$, exists regardless of the number of realization of certain deterministic set of frequencies of harmonic components $\omega_{kq} \in \Omega_{xq}$ in the representation of random processes. Relationship between the sets of realizations $x_{iq}(t) \in X_q \subset \bigcup_q X_q = S_x$ and the system of multiple frequencies $\omega_{kq} \in \Omega_{xq} \subset \bigcup_q \Omega_{xq} = \Omega_x$, generating a basis of subspace establishes $X_q$ the following lemma 1.

It is lemma 1. *Given a finite or countable system of sets, for example,* $S_x = \bigcup_q X_q$ *random stationary in the broad sense of nonergodic processes with discrete spectra* $x_{iq}(t) \in X_q$ *where* $i = 1,2,3,..$ *number of implementation,* $q = 1,2,3,...$ *set number. Each set* $X_q$ *system of sets* $S_x$ *generates a correlation function of the general form* $R_q(\tau) = \sum_{k=1}^{k=\infty} \gamma^2(\omega_{kq})\cos(\omega_{kq}\tau)$ *where* $\omega_{kq} \in \Omega_{xq} \subset \bigcup_q \Omega_{xq}$, $\sum_{k=1}^{k=\infty} \gamma^2(\omega_{kq}) < \infty$ *for everyone* $q$. *If the system sets* $S_x$ *is the union of linearly* *independent sets of random processes, the system sets the angular frequency* $\Omega_x = \bigcup_q \Omega_{xq}$ *generating orthonormal trigonometric basis for each set of random processes* $x_{iq}(t) \in X_q$ *is a semiring.*

It is proof. Rewrite the orthogonality condition for the realization of random processes in complex form

$$M\{\sum_{k=-\infty}^{k=+\infty}\sum_{m=-\infty}^{m=+\infty} C_i(j\omega_{kq})C_i(-j\omega_{md})\exp(j(\omega_{kq}-\omega_{md})t)\}=0, \tag{9}$$

where $\sum_{k=-\infty}^{k=+\infty} C_i(j\omega_{kq})\exp(j\omega_{kq}t)=x_{iq}(t)$ - randomly chosen realization of a random process in a complex form, $\sum_{m=-\infty}^{m=+\infty} C_i(-j\omega_{md})\exp(-j\omega_{md}t)=\bar{x}_{id}(t)$ - other arbitrarily chosen the complex conjugate of realization of a random process $x_{id}(t)$, $q,d=1,2,3,\dots$.

Condition (9) holds for all values of products of coefficients $C_i(j\omega_{kq})\ C_i(-j\omega_{md})$, if for all $q \neq d$ and all $k,m=1,2,3,\dots$ $\omega_{kq}-\omega_{md}\neq 0$. This condition is always satisfied if the system of sets of angular frequency $\Omega_x = \bigcup_{q=1}^{q=l}\Omega_{xq}$ is a semiring, which suppress the pair wise for all $q \neq d$ is empty is empty $\Omega_{xq}\bigcap\Omega_{xd}=\varnothing$ [11. 39]. Lemma 1 is proved.

## 6. Properties of linearly dependent random processes

In the multidimensional control objects on the arbitrarily chosen output synchronously with the input processes $\tilde{x}_{iq}(t)$ there is an output random process

$$\tilde{y}_{ip}(t)=\sum_{q=1}^{q=l} A_{iq}k_{qp}+m_{ip}(t), \tag{10}$$

where the operators $A_{iq}$ are generated by precise linearly independent random components of the input signals $x_{iq}(t)$.

From equation (10) we find randomly selected dedicated channel controls for which the image and preimage associated operator equation $A_{iq}k_{qp}=y_{iqp}(t)$. Property for a set of orthogonal frequency of the harmonic basis for the dedicated control channel is defined by lemma 2.

Lemma 2. *If a linear stationary operator of convolution type* $A_{iq}k_{qp}=y_{iqp}(t)$, *defined in the Hilbert space of almost periodic functions in the sense of Besicovitch* $B_2(-\infty,+\infty)$, *generated by a stationary nonergodic process with a discrete spectrum* $x_{iq}(t)\in X_q \subset B_2(-\infty,+\infty)$, *where* $i=1,2,3,\dots$ *- number of implementation,* $q=1,2,3\dots l$ *- set number, shows the weight function* $k_{qp}(\tau)\in B_2(-\infty,+\infty)$, *where* $p=1,2,3,\dots d$ *on the set of processes* $y_{iqp}\in Y_p \subset B_2(-\infty,+\infty)$, *the n-dimensional countable set of frequencies* $\omega_{cq}\in\Omega_{xq}$, $c=1,2,3,\dots n$, *generated by a countable-dimensional harmonic basis in a separable subspace* $X_q$ *and n-dimensional countable set of frequencies* $\nu_{cp}\in\Xi_{yp}$, $c=1,2,3,\dots n$, *generated by a countable-dimensional harmonic basis in a separable subspace* $y_{iqp}(t)\in Y_p$ *are equal* $\Omega_{xq}=\Xi_{yp}$.

It is proof. The proof is an obvious consequence of the properties of completely continuous operators, the normal form. It is known [12, pp. 186] that in a Hilbert space $B_2(-\infty,+\infty)$, the operator of convolution type $A_{iq}k_{qp} \equiv M\{x_{iq}(t-\tau)k_{qp}(\tau)\} = y_{iqp}(t)$, which is completely continuous normal type operator mapping $B_2(-\infty,+\infty)$ into itself.

By definition, each basis vector is transferred by the operator of a normal form with a coefficient equal to the eigenvalue [17, pp. 203]. It follows that the orthonormal system of basis functions and transform the image of the same.

Known [17, pp. 214] that each eigenspace corresponding nonzero eigenvalue of a completely continuous symmetric operator $A_{iq}$ is finite. Therefore, the set $n$-dimensional angular frequency $\Omega_{xq} = \Xi_{yp}$. Lemma 2 is proved.

## 7. System with independent input actions.

It is theorem 5. *Given mathematical model of a multidimensional linear stationary control object that contains multiple inputs and multiple outputs. Each output is additively associated with each input through the control channel. The control channel is the operator of convolution type, defined in the Hilbert space of almost periodic functions* $B_2(-\infty,+\infty)$ *accurate input random process and the weighting function given by the number input and output number. Each input accurate and precise output processes are distorted by various additive noise and form approximate the original data. The exact components of the process and distort their noise is uncorrelated between a stationary nonergodic random processes with discrete spectra, belonging to the Hilbert space of almost periodic functions in the sense of Besicovitch* $B_2(-\infty,+\infty)$. *If the system sets all of the approximate input signals is linearly independent random processes, then all sets of projections on harmonic bases approximate the input and output processes can be approached to allocate the exact components of the projection of input actions and projection of the exact components of the output processes for each control channel.*

It is proof. Arbitrarily select the input $q$ on which there is an exact realization of a random process, distorted by an additive uncorrelated noise $\tilde{x}_{iq}(t) = x_{iq}(t) + n_{iq}(t)$. According to Lemma 1 the set of frequencies of harmonic components (projections) in the representation of this random process $\tilde{\Omega}_{xq} = \Omega_{xq} \bigcup \Omega_{nq}$ is a semiring $\Omega_{xq} \bigcap \Omega_{nq} = \varnothing$.

In the multidimensional control objects on the arbitrarily chosen output $p$ synchronously with the input processes $\tilde{x}_{iq}(t)$ observed output random process (10). The observed random process $\tilde{y}_{ip} \in \bigcup_{q} Y_{qp} \bigcup M_p$, where $A_{iq}k_{qp} \in Y_{qp}$, $m_{ip} \in M_p$ is an element of a system of linearly independent sets. By Lemmas 1 and 2, the system sets the frequency of the harmonic components (projections) of a random process observed at the output $p$, $\tilde{\Omega}_{yp} = \bigcup_{q=1}^{q=l} \Omega_{xq} \bigcup \Omega_{mp}$ is a semiring $\Omega_{xq} \bigcap \Omega_{mp} = \varnothing$ for all $q$, $\Omega_{xq} \bigcap \Omega_{xd} = \varnothing$ for all $q \neq d$.

We find the source data as the intersection of the sets of frequencies $\tilde{\Omega}_{yp} \bigcap \tilde{\Omega}_{xq}$ approximate the original data, the observed input $q$ and $p$ output,

$$\nu_{kqp} \in \Omega_{qp} = \tilde{\Omega}_{yp} \bigcap \tilde{\Omega}_{xq} = (\bigcup_{c=1}^{c=l} \Omega_{xc} \bigcup \Omega_{mp}) \bigcap (\Omega_{xq} \bigcup \ \Omega_{nq}) = \Omega_{xq} \text{ at } c = q,$$

where $k = 1,2,3,...n$. The resulting intersection of the frequencies coincides with the set frequency, which determines the harmonic basis for accurate component of the input to the $q$-th input. In the limit of infinite observation time, you can find the projection $x_{iq}(j\nu_{kqp})$ $= (\tilde{x}_{iq}(t), \exp(-j\nu_{kqp}t))$ exact spectra of the first implementation of the input $x_{iq}(t)$ at the inlet $q$ and the exact components of the projection of the second implementation $y_{iqp}(j\nu_{kqp})$ $= (\tilde{y}_{ip}(t), \exp(-j\nu_{kqp}t))$ , observed at the output $p$ , the harmonic basis, determined by a set of frequencies $\nu_{kqp} \in \Omega_{xq}$ accurate projections on the background of additive noise $n_{iq}(t)$ and $m_{ip}(t)$. Theorem 5 is proved.

Dedicated accurate projection of input and output components of the processes on the set of input data can be used to address a number of applications in the presence of interference. If we solve the problem of identification of dynamic performances, given the relationship $y_{iqp}(j\nu_{kqp}) = w_{qp}(j\nu_{kqp})x_{iq}(j\nu_{kqp})$, You can find the frequency transfer function of the control channel between the $q$-th entrance and exit to $p$-th $w_{qp}(j\nu_{kqp})$ on the set of frequencies $\nu_{kqp} \in \Omega_{xq}$. If we solve the problem of indirect measurements, such as the acceleration input $q$ against the background of additive noise $n_{iq}(t)$, then filtered from the noise $n_{iq}(t)$ , $i$-i implementation of acceleration is given by $x_{iq}(t) = \sum_{k=-n}^{k=+n} x_{iq}(j\nu_{kqp})\exp(j\nu_{kqp}t)$ .

Apparently, the processes with discrete spectra are divided into periodic, non-periodic and almost periodic processes. If the process is periodic, there exist finite values of time n$T, n = 1,2,3,..$, where we have, for example, process $x(t) = x(t + nT)$. Periodic process is predictable. It suffices to find the minimum value of the final period $T$. The frequency of the first harmonic and other frequency harmonics in the Fourier series are comparable (the frequency of the first harmonic is a measure for the multiple frequencies of the other harmonics).

The main difference between non-periodic processes is that the frequencies (periods), the harmonics are incommensurable. For each harmonic is necessary to determine the angular frequency, which is expressed in non-periodic infinite decimal. Instrumental and methodological errors limit the accuracy of the frequencies. Find the exact mathematical model of the observed non-periodic process is impossible. Therefore, it is unpredictable throughout the time axis and on this basis is a model of a random process.

Almost periodic processes occupy an intermediate value. Frequency harmonic components can be set arbitrarily. These may include harmonic components with multiple and disparate frequencies, ie contain periodic and aperiodic components. Due to the non-periodic component of almost periodic processes are unpredictable, as well as non-periodic processes.

There are various methods for determining the approximate frequency of the harmonic orthogonal basis of the observed random processes, such as input signals $\tilde{\Omega}_{xq} = \Omega_{xq} \bigcup \Omega_{nq}$. Among them are well known method of determining the amplitude of the current spectrum [18, pp. 111]. Seeking an array of frequencies is determined by the frequency at which the current amplitude spectrum has a local maximum.

## 8. System with correlated input actions

Described the frequency control method of noise is applied only in cases where all input actions are linearly independent. Moreover, for all inputs must be executed Conditions $\tilde{\Omega}_{xq} \bigcap \tilde{\Omega}_{xd} = \varnothing$ at $xq \neq xd$ , $q, d = 1,2,3,..\, l$. If at least one pair of different inputs of this condition is not satisfied, then it indicates that the test object is an object of control with correlated input actions.

Degree of connection between the exact components of the input signals in real form, observed on different inputs, for example, $x_q$ and $x_d$ can be estimated by the cosine of the angle between them $\cos(x_q, x_d)$ [11, pp. 135]. In this case, there are three cases. If $|\cos(x_q, x_d)| = 1$, the signals observed at the entrances $x_q$ and $x_d$, coincide up to sign. If after determining the proper dimension of the object controls the number of entries it appears that all remaining inputs, the condition $\cos(x_q, x_d) = 0$ , this indicates that the condition of theorem 5 holds.

In real terms, as a rule, for all inputs or a subset of inputs, the condition $|\cos(x_q, x_d)| < 1$ | for $q \neq d$ , which indicates the correlation of the input signals. In this case, the above method of dealing with noise is not applicable. The described method of interference mitigation can be applied if, instead of the original system with correlated effects using the subsystem with uncorrelated input actions. Based on theorems 6 a subsystem can be obtained.

It is theorem 6. *If the system is real accurate component of the input signal is a system that satisfies* $|cos(x_q, x_d)| < 1$ *when* $q \neq d$ *for all inputs or a subset of input signals, and for another subset of the input signals* $cos(x_q, x_d) = 0$, *then the original system input signals that are distorted correlated and uncorrelated noise, we can distinguish linearly independent subsystem of input signals, distorted only by additive uncorrelated noise.*

It is proof. Change the previously used model inputs so that would condition$|\cos(x_q, x_d)| < 1$ holds for all $x_q \neq x_d$ . Such a system of sets of input actions denoted as

$$\widehat{S}_{xnf} = \bigcup_{q=1}^{q=l} \widehat{X}_q, \text{ где } \widehat{x}_{iq}(t) \in \widehat{X}_q = \{\widehat{x}_{iq}(t) : \widehat{x}_{iq}(t) = x_{iq}(t) + n_{iq}(t) + \sum_{c=1}^{c=l} f_{iqc}(t),\ c \neq q\},$$

$x_{iq}(t)$ - the exact components of the input signal and the random additive noise $n_{iq}(t) \in N_q$ at the input of q,

$f_{iqp}(t)$ - $i$ - th realization of the random correlated noise in the form of a stationary random function of communication, which is also observed at different entrances $q$ and $p$ .

All the components of each realization of random processes are linearly independent (uncorrelated) with each other. Therefore, only the stationary random processes $f_{iqp}(t)$ that occur at different inputs $q$ and $p$ generate non-trivial cross-correlation function $E\{\widehat{\widehat{x}}_{iq}(t-\tau)\widehat{x}_{ip}(t)\} = R_{fqp}(\tau)$ .

The set of frequencies of harmonic components of random processes at each input $\widehat{\Omega}_{xq} = \Omega_{xq} \bigcup \Omega_{nq} \bigcup_{c=1}^{c=l} \Omega_{fqc}$ , $q = 1,2,3,...l$, for all $c \neq q$. The set of frequencies of harmonic components of the function of communication $\Omega_{fqp}$ , observed between the inputs $p$ and $q$ are by definition the intersection of sets $\widehat{\Omega}_{xp} \bigcap \widehat{\Omega}_{xq} = \Omega_{fqp}$ .

Obtain the system sets the corner frequency of the harmonic components of all the functions of communication, $F = \bigcup \Omega_{fqp}\ q, p = 1,2,3,...l$, for all $p \neq q$ . For an arbitrarily selected input $q$ will find the difference [11, c. 13] systems of sets of angular frequency

$$\widehat{\Omega}_{xq} \setminus F = \Omega_{xq} \bigcup \Omega_{nq} \bigcup_{c=1}^{c=l} \Omega_{fqc} \setminus F = \Omega_{xq} \bigcup \Omega_{nq} = \widetilde{\Omega}_{xnq} \tag{11}$$

Performing the operation (11) for all $q = 1,2,3,...l$ obtain a system of sets of frequencies $\bigcup_q \widetilde{\Omega}_{xnq}$ which is a semicircle, as $\widetilde{\Omega}_{xnp} \bigcap \widetilde{\Omega}_{xnq} = \varnothing$ for all $p \neq q$ . Thus, applying the subtraction (11) to the original system input correlated effects $\widehat{S}_{xnf}$ You can select from it a different system $\widetilde{S}_{xn}$ with uncorrelated effects, which is determined by the system sets the frequency of all harmonics $\tilde{\Omega}_{xn} = \bigcup_q (\Omega_{xq} \bigcup \Omega_{nq})$ exact components of linearly independent input signals $\bigcup_q \Omega_{xq}$ and all harmonic frequencies of additive noise $\bigcup_q \Omega_{nq}$ . Theorem 6 is proved.

On an arbitrarily selected output $p$ observed $i$ - i realization of a random process generated by the original correlated system of input actions

$$\hat{y}_{ip}(t) = y_{ixp}(t) + \vartheta_{ip}(t) + m_{ip}(t), \quad (12)$$

where $y_{ixp}(t) = \sum_{c=1}^{c=l} x_{ic} * k_{cp}$ - the exact component of the output process, generated by independent input actions, $x_{ic}(t)$, $x_{ic} * k_{cp}$ - operator convolution type defined in a Hilbert space $B_2(-\infty,+\infty)$, $k_{cp}$ - weighting function of the control channel between the input c and output $p$;

$\vartheta_{ip}(t) = \sum_{r=1}^{r=l-1} \sum_{c=r+1}^{c=l} f_{irc} * (k_{rp} + k_{cp})$ - component of the signal generated at the output $p$ function of communication $f_{irc}$ which is observed between the inputs $r$ and $c$, $f_{irc} * (k_{rp} + k_{cp})$ - the symbol of convolution type defined in the Hilbert space $B_2(-\infty,+\infty)$, $k_{rp}$, $k_{cp}$ - weighting function of control channels between, respectively, inputs $r, c$ and output $p$ and $(y_{ixp}, \vec{\vartheta}_{ip}) = 0$, as $x_{ic}$ and $f_{irc}$ - are linearly independent random processes for all $c$ and $r$, $m_{ip}(t)$ - $i$ -i implementation of the additive noise, which satisfies the conditions of independence $(y_{ixp}, \overline{m}_{ip}) = (\vartheta_{ip}, \overline{m}_{ip}) = 0$.

Realizations of random processes at different outlets, for example, $p$ and $d$ processes are correlated due $(\hat{y}_{ip}, \hat{\bar{y}}_{id}) \neq 0$ to mutual coupling between inputs. However, the random process observed at the output $p$ is the sum of three mutually orthogonal processes (12). By lemma 1, set the frequency of the harmonic components of the process $\hat{y}_{ip}(t)$ is a semiring of sets of frequencies $\hat{\Omega}_{yp} = \bigcup_{c=1}^{c=l} \Omega_{xc} \bigcup \Omega_{\vartheta p} \bigcup \Omega_{mp}$ . With the set of frequencies $\tilde{\Omega}_{xnq}$, defined by the formula (11), we find the set intersection frequency

$\nu_{kqp} \in \hat{\Omega}_{yp} \bigcap \tilde{\Omega}_{xnq} = (\bigcup_{c=1}^{c=l} \Omega_{xc} \bigcup \Omega_{\vartheta p} \bigcup \Omega_{mp}) \bigcap (\Omega_{xq} \bigcup \Omega_{nq}) = \Omega_{xq}$ for an arbitrarily chosen pair of indices of inputs $q$ and outputs $p$.

Thus, by the operation of the selection of independent components (11) input signals $\tilde{S}_{xn}$ from the original correlated system of input actions $\hat{S}_{xnv}$ frequency control method of noise in systems with independent input actions can be generalized to cases where input actions are correlated effects.

## 9. An example of solving the problem of identification

The proposed method of dealing with noise is used in solving the problems of identification of dynamic characteristics of the multidimensional control objects in a class of linear stationary models. In particular, the initial data, observed in the regime of one automatic landing aircraft (class IL -96), the mathematical model of a control channel in a class of linear stationary models between input $x_1$(t) - a given pitch, at the input-wire control system for pitch, and output $y(t)$ - the actual angle pitch plane. On the pitch is affected by other inputs: the position of thrust lever $x_2$ (t), the angular position of the flaps $x_3(t)$ and angular position of the stabilizer $x_4$(t).

Synchronous records listed the output and input signals were provided in the form of 274 discrete samples in 0,5 sec. (ZAK_51_IL_V30). Inputs $x_3(t)$ and $x_4$(t) signals are relay type, derivatives of which are described by finite functions on the whole line in the form of rectangular pulses. As discussed above (Section 4), the proposed method is applied to random processes, which have a continuous second derivative. For signal relay type, this condition is not satisfied, making it difficult for the solution of problem identification through $x_3(t)$ - $y(t)$ and $x_4(t)$ - $y(t)$.

Known method for solving the problem of identification of the dynamics of multidimensional control systems reduces to the solution of integral equations of the first kind [19, pp. 22]. System of equations can be obtained only for the control objects for which the number of inputs equals the number of outputs. Baseline data are presented in the form of autocorrelation functions of all input signals, cross correlation functions between all the different inputs and cross correlation functions between all input and all output signals. Observed signals are stationary random fragments in the broad sense of processes that do not possess the ergodic property. This fact excludes the receipt of the correlation functions by averaging over time of individual implementations. Obtaining the correlation functions over the set of realizations is hampered for the reason that it is impossible to perform many aircraft landings in the same conditions, such as weather. In addition, all aircraft like objects have a dimension of control over all inputs greater than that of all the exits. For these reasons, the task set can not be solved by a known method.

The task of identifying the dynamic characteristics of the channel management has been solved by a new method proposed in the following sequence. First by local maxima of the current amplitude spectra of the observed signals were determined frequency harmonic components simultaneously observed signals, which are secondary source data. Received the original system sets the frequency of the harmonic components of input and output signals $\cup_i \Omega_{xi}, \Omega_y$, where $\Omega_{xi}$ - set the frequency of the harmonic components of input signals on the $i$ -th entry, $i = 1,2,3,4$, $\Omega_y$ - set the frequency of the harmonic components of output signal $y(t)$.

Theorem 5 implies that if the system sets the approximate data observed at all entry points is linearly independent random processes, then all sets of projections on harmonic bases approximate the input and output processes close, you can select the exact components of the projection of input actions and projection of the exact components of output processes. Typically, the system sets the input signal does not satisfy Theorem 5. Therefore, by subtraction of sets (11) of the resulting system sets the frequency of the harmonic components of input signals is a subsystem of the sets of frequencies only independent of the harmonic components of input signals, which are tertiary source data. Software, this operation is performed as follows. Chosen criterion under which two frequencies are coincident. If $|\omega_{xkp} - \omega_{xmq}| \le \delta$ , where $\delta = 2\pi/T$, $T$ -duration of implementation, for $\forall(xk \neq xm)$, $q \in [1, r_m]$, $p\epsilon[1, r_k]$, $p\epsilon[1, r_k]$, $r_k -$ dimension array of frequencies $\Omega_{xk}$, $r_m$ - the dimension of the array of frequencies $\cup_i \overline{\Omega}_{xi}$, $k, m \in [1,4]$ then the frequencies are coincident and are removed from each source array frequency input signals. As a result, only selected frequencies of arrays of independent harmonic components of input signals $\cup_i \overline{\Omega}_{xi}$, as a rule, c smaller dimension. For example, the original dimension of the array of frequencies at the first entrance was $r_1 = 137$, and the dimension of the array of frequencies independent of harmonic components on the same first entry was $\bar{r}_1 = 32$.

By lemma 2, for a given control channel $x_1 - y$ is defined by a set of frequencies of harmonic components that generate and describe the forced motion in the channel controls on the criterion that the frequencies of harmonic components of input and output signals of the selected channel management $\Omega_{x1y} = \{\omega_\rho = \omega_i \in \Omega_{x1y}: |\omega_i - \omega_k| \le \delta, \forall \omega_i \in \overline{\Omega}_{x1}, \forall \omega_k \in \Omega_y, i \in [1, \bar{r}_{x1}], k \in [1, r_y]\}$, where $r_y = 94$ - the dimension of the resulting array of frequencies harmonic output $\Omega_y, \rho \in [1, d]$ , where $d$ - the dimension of the matching frequency harmonic components, generators, and describing the forced motion in the channel $x_1 - y$. In this example, it was found that $d = 13$. A subset of frequencies are $\Omega_y$ quaternary source data.

For all received frequencies $\omega_\rho \in \Omega_{x1y}$ determined by the Fourier exponents of the input signal $S_{x1}(j\omega_\rho) = a(\omega_\rho) + jb(\omega_\rho)$ and the output signal $S_y\left(j\omega_\rho\right) = \gamma(\omega_\rho) + j\beta(\omega_\rho)$. At the lowest frequency $\omega_{\rho 1}$ sets of frequencies $\omega_{\rho 1} \in \Omega_{x1y}$ astatizma determines the order in the system for the control action. For what is the value of the frequency transfer function at the lowest frequency $W(j\omega_{\rho 1}) = S_y(\mathrm{j}\omega_{\rho 1})/S_x(\mathrm{j}\omega_{\rho 1}) = P + jQ$ . Depending on the position of the point $W(j\omega_{\rho 1})$ in the complex plane is determined by the order astatizma on the control action. Usually, the specified point can be located in the first, second or third quadrant of the complex plane, which corresponds to the zero, first or second order astatizma $p_a$ [20, s.300].

Obtained by Fourier exponents and order astatizma (in this example $p_a = 1$) can make a number of systems of algebraic equations on the 2 nd to the $d$ - th order, which are the Fourier transforms of systems of linear ordinary differential equations with unknown constant coefficients of the 2 nd to $d$ - order. For example, a system of algebraic equations, given the current $q$ -th order, where $2 \le q \le d$, is represented as

$$\sum_{k=0}^{k=q} (j\omega_{p1})^{k+p_a} T_{k+p_a} S_y(j\omega_{p1}) = S_{x1}(j\omega_{p1}),$$

$$\dots\dots\dots\dots\dots\dots\dots\dots\dots\dots\dots\dots\dots\dots\dots$$

$$\sum_{k=0}^{k=q} (j\omega_{pc})^{k+p_a} T_{k+p_a} S_y(j\omega_{pc}) = S_{x1}(j\omega_{pc}),$$

where c$\in [1, q]$ .

From this system of equations should be two systems of equations for real and imaginary components. Solving these equations successively from 2 nd to $d$-th order with respect to coefficients $T_{k+p_a}$, we find that the maximum order of the system of equations $q = q_{max}$, at which the system of equations is consistent. To solve this problem it turned out that the dynamic characteristics of the selected channel management, taking into account the elastic strain aircraft structure can be described as a first approximation of ordinary differential equations 9-th order with constant coefficients: $T_0$= 0; $T_1$ =- 2,8531; $T_2$ =- 14,5651; $T_3$=- 5,4709; $T_4$ =- 3,9165; $T_5$ =- 1,0721; $T_6$=- 0,2584; $T_7$ =- 0,06054; $T_8$ =- 0.004845; $T_9$ =- 0,001023. Negative signs of the coefficients of the differential equation due to the rule of signs taken in the aerodynamics, to the tilt direction of the elevator and the longitudinal axis of the aircraft in pitch [21, pp. 6,. 25].

In a similar way was obtained a differential equation for the other control channel, which determines the effect of engine power control on the pitch. More information about a new method for identification of dynamic characteristics is published on http://asvt51.narod.ru.

The obtained differential equations in two-channel pitch control allows the system to find the control of the aircraft taking into account the elastic deformation the airframe in the space of states.

## 10. Conclusion

1. In the modern theory of random wide-sense stationary processes are the fundamental results of empirical research on E. Slutsky and the results of theoretical studies of A. Khintchine and Wiener respect to the representation of correlation functions of random stationary processes. Representation of themselves random stationary processes in a broad sense have been derived from empirical studies E. Slutsky and as a result of theoretical research in the papers of Kolmogorov, Cramer and other scientists. In theoretical studies suggest that, firstly, the independent variable is the angular frequency is a continuous variable. Secondly, the probability distribution function of angular frequency includes a step and continuous components. It turned out that the results of empirical studies contradict the theoretical studies. It was conjectured that the contradictions are the result of incorrect assumptions. It turned out that theoretical studies were performed at incorrect assumptions. Following the theoretical studies with reasonable assumptions, it turns out that random stationary in the broad sense has a unique representation as a Fourier series, defined in the Hilbert space of almost periodic functions that have only a discrete spectrum, has the ergodic property of the first order and has no ergodic property of the second order. In this empirical and theoretical studies do not contradict each other.

2. To combat the uncorrelated noise in Wiener instead of primary source data - observable input and output random processes using secondary source data in the form of deterministic correlation functions of random processes. Due to the non-ergodicity of random processes, second order, this method of dealing with noise is limited in the application of those rare cases where there is a possibility of correlation functions by averaging over multiple realizations, and the set of processes generated by the forced motion in the channel of the control object is not distorted by even small interference. In this regard, has been tasked to develop a method for interference mitigation in circumstances where the original data given in the form of single realizations of random nonergodic processes.

3. A new method for interference mitigation, based on the separation of double frequency harmonic components of the sets of bases observed signals of multidimensional control systems. As a secondary deterministic data are encouraged to use the deterministic set of harmonic frequencies of bases observed at the inputs and outputs. In many cases the set of input signals is correlated with each other. This leads to a correlated noise output signals with input signals. To overcome this limitation is proposed tertiary treatment of the original input data. Harmonic frequencies of the sets of input signals are removed coincident frequencies. The task of identification of dynamic characteristics for dependent input signals is reduced to a simpler problem of identifying the dynamic characteristics of the independent input signals. To identify the selected control channel used tertiary data input signals and the secondary data of output signals. The quaternary initial data there is a subset of conterminous frequencies between secondary and tertiary initial source data. Quaternary source data contains only a subset of the frequencies of harmonic components of accurate baseline data. The identification problem and the problem of filtering the input and output signals from noise is solved for the selected control channel only on a subset of the frequencies of Quaternary data.

4. The described method of dealing with noise is used for solving the problems of identification of dynamic characteristics of the control channels of multidimensional control objects in a class of linear stationary models. Is an example of solving the problem of identification of the control channel - given pitch angle, at the input-wire control system pitch - the actual pitch of the aircraft in the class of linear stationary models based on initial data obtained during one automatic landing (class of IL-96). Known method of solving the problem proved to be ineffective for two reasons. Firstly, in the test facility management is not satisfied the restriction on the number of inputs equal to the number of outputs. At one output simultaneously affect signals observed in the four entrances. This eliminates the need to obtain a system of integral equations of the first kind. Secondly, it is impossible to obtain the correlation functions by averaging the random signals to multiple implementations. The proposed method of dealing with noise in the conditions described above allowed us to obtain dynamic response of the control channel in a first approximation, taking into account the elastic strain aircraft structure in the form of ordinary differential equations 9-th order with constant coefficients. The control channel is first order astatizma. After that, it was decided to identify the dynamic characteristics of the other control channel in the coordinates - the position of control levers engines - the pitch plane. The dynamics of this control channel is described as a first approximation, taking into account the effect of elastic deformation the airframe ordinary differential equation of fifth order with constant coefficients.

The information obtained on the differential equations for the two channels control the pitch control system allows you to build a pitch in the space of states.